\documentclass[12pt]{article}
\usepackage{amsmath}
\usepackage{amsthm}
\newtheorem{theorem}{Theorem}
\newtheorem{corollary}{Corollary}

\newtheorem{proposition}{Proposition}
\newtheorem{lemma}{Lemma}
\newtheorem{conjecture}{Conjecture}

\title{On Thin Sets of Primes Expressible as Sumsets} 
\author{Ernest S. Croot III and Christian Elsholtz}
\date{\today}
\begin{document}
\maketitle

\section{Introduction.}

In this paper we will use the following notation.  Given a set of positive 
integers $S$, we let $S(x)$ denote the number of elements in $S$ that are 
$\leq x$, and we let $|S|$ denote the total number of elements of $S$. 
Given two sets of positive integers $A$ and $B$, we denote the
{\it sumset} $\{ a+b\ :\ a \in A, b\in B\}$ by $A+B$; and so, the
number of elements in $A+B$ that are $\leq x$ will be $(A+B)(x)$.
For a finite set of integers $J$, and integers $q \geq 2$ and
$r$, we let $J(r,q)$ denote the set of elements of $J$
which are $\equiv r \pmod{q}$. 
We will also use Vinogradov's notation:  The statements ``$f(x) \ll g(x)$''
and ``$g(x) \gg f(x)$'' are both equivalent to ``$f(x) = O(g(x))$'';
and, we will use ``$f(x) \ll_y g(x)$'' to indicate 
that the implied constant in the big-oh depends on a parameter $y$.
Finally, 
by the statement $f(x) \sim g(x)$ we mean that 
$$
\lim_{x \to \infty} {f(x) \over g(x)}\ =\ 1.
$$

An old conjecture of Ostmann, which is sometimes called the 
`Inverse Goldbach Problem', asks whether there is an additive
decomposition of the primes,
with at most finitely many 
exceptions (see \cite{Ostmann}, p.~13); that is, do there exist 
sets of positive integers $A$ and $B$, each with at least two elements,
such that 
$$
{\rm for\ } n > x_0\ {\rm (for\ some\ }x_0{\rm )},\ \ n \in A+B\ \iff\ n{\rm\ is\ 
prime}.
$$ 

Even though this question has withstood attack by several
mathematicians, there has been much recent progress.  For instance, 
Wirsing in \cite{Wirsing}, 
Pomerance, S\'{a}rk\"{o}zy and Stewart in 
\cite{PomeranceandSarkozyandStewart}, 
Hofmann and Wolke in \cite{Hofmann} and Bshouty and Bshouty in 
\cite{BshoutyandBshouty}
have shown that if 
such a decomposition exists, then one has the following estimates
on the counting functions
$$
{x \over \log x} \ll A(x) B(x) \ll x,
$$
Elsholtz has shown in \cite{Elsholtz2} that 
$$
{\sqrt{x} \over \log^5x}\ \ll\ \min(A(x), B(x))\ \leq\ \max(A(x),B(x))
\ \ll\ \sqrt{x}\log^4x,
$$
which improves on another result of Hofmann and Wolke from
\cite{Hofmann}.  Furthermore, Elsholtz uses these results to 
show that $B$ (or $A$) cannot be the 
sum of two other sets, each with at least two elements; that is,
Elsholtz has solved a ternary analogue of the 
above conjecture of Ostmann.  

It seems conceivable that sieve methods alone will not solve the Ostmann
problem, but that some additional insight into the structure of sumsets
is needed.  This paper is a step towards this.
\bigskip

In his proof of the above mentioned result, 
Elsholtz makes strong use of the fact that 
$$
(A + B + C)(x)\ \gg\ {x \over \log x}.
$$
This leads one to wonder whether this constraint can be weakened
somewhat.  To be more specific:  
\begin{quote} Does there exist $\kappa > 1$ and sets
of positive integers $A,B,C$, each with at least two elements, 
such that $A+B+C$ is a set of primes with 
\begin{equation}\label{hypothesis_inequality}
(A + B + C)(x)\ \gg\ {x \over \log^\kappa x}? 
\end{equation}
\end{quote}
And the answer to this question is:  Yes.  The Hardy-Littlewood conjecture can
be used to give solutions for every $\kappa \geq 3$.  
Before we show how, we give
here the form of the conjecture we will need (see \cite{HardyandLittlewood}):
\bigskip

\noindent {\bf Hardy-Littlewood Conjecture:}  Suppose that 
$a_1 < a_2 < \cdots < a_k$ is a sequence of integers such that the polynomial
$(x+a_1)(x+a_2)\cdots (x+a_k)$ has no fixed prime divisors.  Then,
$$
\#\{ n \leq x\ :\ n+a_1,n+a_2,...,n+a_k\ {\rm are\ all\ prime} \}\ \sim\ 
C(a_1,...,a_k){x \over \log^k x},
$$
where $C(a_1,...,a_k)$ is some constant which depends only on $a_1,...,a_k$.
\bigskip

\noindent  Now, suppose that $A = B = \{1,7\}$, and let $C$ be the set of all positive 
integers $n$ such that $n+2, n+8, n+14$ are all prime.  
Then, $A+B+C$ consists 
entirely of primes, and since $(x+2)(x+8)(x+14)$ 
has no fixed prime divisors, assuming the Hardy-Littlewood conjecture we get that
$$
C(x)\ \sim\ C(2,8,14) {x \over \log^3 x};
$$   
and so,
$$
(A+B+C)(x)\ \gg\ {x \over \log^3 x},
$$
which means that our question above has an affirmitive answer for all 
$\kappa \geq 3$.

So, if we are to have any hope of extending Elsholtz's work to show that there
are no triples $A,B,C$ where $A+B+C$ contains {\it many} primes 
(but not almost all primes), we have to account for the above ``obstruction'' 
arising from the Hardy-Littlewood conjecture.  The following are all the cases
where one can apply the Hardy-Littlewood conjecture to construct sets $A,B,C$
such that (\ref{hypothesis_inequality}) holds:
\begin{equation}\label{pair_inequalities}
{\rm Either\ }|A+B| \leq \kappa,\ {\rm or\ } |A+C| \leq \kappa,\ {\rm or\ } 
|B+C| \leq \kappa.
\end{equation}
This now leads us to the following general conjecture:

\begin{conjecture}\label{ternary_conjecture}  
Suppose $A,B,C$ are sets of positive integers with at
least two elements each, such that
$A+B+C$ consists entirely of primes.
If 
$$
(A + B + C)(x)\ \gg_\kappa\ {x \over \log^\kappa x},
$$
then (\ref{pair_inequalities}) holds.
\end{conjecture}

In this paper we do not quite prove this conjecture, although we believe
that it is true.  One additional, technical assumption about the sets $A,B,C$
is needed for our proof; basically, we need that there are not ``too many'' 
primes $p$ which have ``too many'' solutions $p = a+b+c$, $a \in A, b\in B, c \in C$.
To state this technical assumption, we need the following definition:  
\bigskip

\noindent {\bf Definition.  }  For a given collection of sets 
$A_1,A_2,...,A_k$ let $r(n;A_1,...,A_k)$ denote the number of solutions
to
$$
n = a_1 + \cdots + a_k,\ a_i \in A_i\ {\rm for\ all\ }i=1,2,...,k.
$$
We say that the collection of sets $A_1,...,A_k$ is {\it regular} if
and only if for every $\epsilon > 0$, there exists $D > 0$ such that
for $x$ sufficiently large,
\begin{equation}\label{regular_inequality}
\sum_{n \in S} r(n;A_1,...,A_k)\ <\ 
\epsilon A_1(x)\cdots A_k(x),
\end{equation}
where 
\begin{equation}\label{Sx_regular}
S\ =\ \{ n \in A_1 + \cdots + A_k\ :\ n\leq x,\  
r(n; A_1,...,A_k) > \log^D x\}.
\end{equation}
\bigskip

One easy consequence of the fact that $A_1,...,A_k$ is regular is
the following:
\begin{lemma}\label{regular_lemma}  Suppose $A_1,...,A_k$ is regular.  
Then, there exists a constant $E > 0$ such that for $x$ sufficiently
large,
$$
A_1(x) A_2(x) \cdots A_k(x)\ <\ (\log^Ex)\ (A_1 + \cdots + A_k)(x).
$$
\end{lemma}
\noindent A proof of this lemma can be found in Section \ref{technical_section}.
\bigskip

Our Main Theorem is as follows:
\bigskip

\begin{theorem} {\bf (Main Theorem)} \label{main_theorem}
Conjecture \ref{ternary_conjecture} holds if we assume $A,B,C$ is a  
regular triple; that is, if $A,B,C$ is a regular triple of sets 
of positive integers such that $A+B+C$ is a set of primes, 
and
\begin{equation}\label{main_theorem_hypothesis}
|A|,\ |B|,\ |C| \geq 2,\ {\rm and\ }
(A + B + C)(x)\ \gg\ {x \over \log^\kappa x},
\end{equation}
then (\ref{pair_inequalities}) holds.
\end{theorem}
\bigskip

\noindent {\bf Note:} The conclusion of this theorem can {\it possibly} be
proved under a weaker notion of {\it regularity}:  One can maybe
replace the ``$\log^Dx$'' in (\ref{Sx_regular}) with  

\noindent ``$\exp( \log^{1-o(1)} x)$'', and still have the theorem hold. 
This would require substantial modifications of many parts of the argument,
including Propositions \ref{regularity_prop} and \ref{few_solutions_prop},
and Corollary \ref{few_solutions_corollary}.

Perhaps the most interesting feature of the theorem is 
the many different ingredients which are used to prove it (it looks like
a problem tailor-made for a single application of some sieve method),
and the way they all fit together.  These include:  the Large Sieve,
Brun's Sive, Gallagher's Sieve, the ``probabilistic method'' and
regularity principles (which are used to prove Proposition 
\ref{regularity_prop}), translation invariant principles (which appear
in Proposition \ref{few_solutions_prop} and Lemma 
\ref{translation_invariant_lemma}), and certain ``maximality'' or 
``local-global'' principles (which appear in the proof of Lemma 
\ref{cauchy_davenport_consequence}).

The basic idea of the proof (of Theorem \ref{main_theorem}) is as follows:
We will prove the contrapositive of the Theorem by first assuming that
(\ref{pair_inequalities}) fails to hold. 
Through a combination of Propositions \ref{elsholtz_prop} and
\ref{extract_B_C} (which appear in the next section of the paper) and
some basic combinatorial arguments, we will find subsets of $A \cap [1,x]$, 
$B \cap [1,x]$ and $C \cap [1,x]$, call these subsets 
$\hat A$, $\hat B$ and $\hat C$, which will have certain usable properties.
At this point, the proof will break down into two cases, with case 1 
being where $\min(|\hat A|, |\hat B|) > \kappa$ and case 2 where
this min is $\leq \kappa$.  The most difficult and important case will be 
case 2; and for this case, we will construct subsets of $\hat A,\hat B$
and $\hat C$, call these subsets $S, L^*$ and $C^*$, such that the
following inequalities hold:
\begin{equation}\label{vague_LC_bounds}
|S| \leq \kappa,\ \ 
|S + L^* + C^*|\ \geq\ {|L^*|\ |C^*| \over 2}\ >\ {A(x)B(x)C(x) \over 
\log^E x},
\end{equation}
for some $E > 0$, and
$$
{\sqrt{x} \over \log^{6+\kappa} x}\ \ll\ |L^*|,\ |C^*|\ \ll\ 
\sqrt{x} \log^6x.
$$

Then, we will show that most triples $(a,b,c) \in S \times L^* \times C^*$ 
have the property that, for any integer $r \geq 1$ and some integer $k$ 
(and $x$ sufficiently large), the numbers
$$
a+b+c+k,\ \ a+b+c+2k,\ \ ...,\ \ a+b+c+rk
$$
have very few prime divisors in certain ``long intervals''.  This 
result will follow by showing that the sets $L^*$ and $C^*$ are
approximately ``locally translation invariant'', meaning
that for `many' primes $p \leq \sqrt{x}$, the residue classes
modulo $p$ occupied by $L^*$ will be almost exactly the same as those
occupied by $L^*+k$, $L^* + 2k$, $....,$ and $L^* + rk$; the same
will hold for $C^*$.  The method used to prove this will involve combining
very precise arithmetic information about the sets $L^*$ and $C^*$
together with a variant of Gallagher's Larger Sieve, and will be 
the subject of Lemma \ref{translation_invariant_lemma} within the proof of
Proposition \ref{few_solutions_prop}. 

Next, using Brun's upper bound
sieve, we will show that the number of integers $n \leq x$ where 
$n+k,n+k,...,n+rk$ all have such few 
prime divisors in these ``long intervals'' is $\ll x \log^{-r/2} x$; and so,
from this and (\ref{vague_LC_bounds}) we will deduce
$$
(A+B+C)(x)\ \ll\ |S + L^* + C^*|\log^E x\ \ll\ x \log^{-r/2 + E}x
\ =\ o(x\log^{-\kappa}x),
$$
for $r > 2(\kappa + E)$.  This will contradict the hypothesis of Theorem
\ref{main_theorem}, and so the Theorem will follow. 
\bigskip

Although we indicated earlier how the Hardy-Littlewood conjecture can
be used to produce sets $A,B,C$ which satisfy the hypotheses and conclusion
of this Theorem for $\kappa \geq 3$, we can give a weaker, unconditional result.
Besides the sharpness of the inequalities obtained through the Hardy-Littlewood
conjecture, this result is also weaker in that it only holds for a fixed $x$.

\begin{theorem}{\label{unconditionallowerbound}}
Given integers $1 < \kappa_1 < \kappa_2$, for all sufficiently large $x$ 
there exist sets of positive integers $A,B,C \subseteq \{1,2,...,x\}$, 
with
$$
|A| = \kappa_1,\ |B|=\kappa_2,\ |C| \geq |B| \geq |A|,
$$
such that $A+B+C$ consists entirely of primes, and
$$
|A+B+C| > c_{\kappa_1,\kappa_2}{x \over \log^{\kappa_1\kappa_2} x},
$$
where $c_{\kappa_1,\kappa_2}$ is some constant depending only on 
$\kappa_1$ and $\kappa_2$.
\end{theorem}

\section{Proof of the Main Theorem (Theorem \ref{main_theorem}).}

We will prove the contrapositive of this theorem.
So, let us suppose that (\ref{pair_inequalities}) fails to
hold; that is,
$$
|A+B|,\ |A+C|,\ |B+C|\ >\ \kappa.
$$

In our proof, we will first require a result that is a slight
generalization of a result of C. Elsholtz \cite{Elsholtz2},
as well as a result which allows us to extract subsets  
$$
\hat A \subseteq A \cap [1,x],\ \hat B \subseteq B\cap [1,x],\ {\rm and}\  
\hat C \subseteq C \cap [1,x],
$$
such that $\hat A + \hat B + \hat C\subset (\sqrt{x},2x)$.
These first two Propositions are as follows:
\bigskip

\begin{proposition} \label{elsholtz_prop}  
If $F, G \subseteq \{1,2,...,x\}$, with $1 \leq \delta < |F| \leq |G|$, 
such that $F+G$ consists entirely of primes in $(\sqrt{x},2x)$, and if
$|F||G| \gg x/\log^\delta x$, then 
$$
{\sqrt{x} \over \log^{\delta+6} x}\ \ll\ |F|\leq |G|\ \ll\ 
{\sqrt{x}\log^6x}.
$$ 
\end{proposition}
\noindent {\bf Note:}  The constant $6$ can certainly be improved here,
but such an improvement does not much affect the quality of the main result
in this paper.\bigskip

\begin{proposition}\label{extract_B_C}
If $A,B,C$ is a regular triple, then there exist subsets 
$\hat A \subseteq A \cap [1,x]$, 
$\hat B \subseteq B \cap [1,x]$, 
$\hat C \subseteq C \cap [1,x]$, such that 
\begin{equation}\label{sim_ABC}
|\hat A| \sim A(x),\ |\hat B| \sim B(x),\ {\rm and\ } 
|\hat C| \sim C(x),
\end{equation}
where
$$
|\hat A + \hat B + \hat C|\ \sim\ (A+B+C)(x) 
\gg {x \over \log^\kappa x},
$$
and $\hat A + \hat B + \hat C \subset (\sqrt{x},\infty)$.
\end{proposition}

For $x$ sufficiently large, we may assume that 
\begin{equation}\label{hat_ABC_inequalities}
|\hat A|,\ |\hat B|,\ |\hat C|\ \geq\ 2,\ {\rm and\ }\ 
|\hat A + \hat B|,\  |\hat A + \hat B|,\ |\hat B + \hat C|\ >\ \kappa.
\end{equation}
The first inequality holds since 
$$
|\hat A| \sim A(x) \geq 2,\ |\hat B| \sim B(x) \geq 2,\ {\rm and\ }
|\hat C| \sim C(x) \geq 2;
$$ 
and the second inequaltiy holds for similar reasons. 

For a given $x$, suppose that, without loss of generality,
$$
A(x),\ B(x)\ \leq\ C(x).
$$
Consider the two sets $\hat A + \hat B$ and $\hat C$, and let 
$F$ be the set with the smaller number of elements, and $G$ be the
set with the larger number of elements.  We will show that these two
sets $F$ and $G$ satisfy the hypotheses of Proposition \ref{elsholtz_prop},
and we will use the conclusion of this proposition to show that 
$|\hat A + \hat B|$ is ``large'', which will be important in later arguments.

We first claim that $|G| \geq |F| > \kappa$ for $x$ sufficiently large:  
To see this, we note that
\begin{eqnarray}
|\hat C|^3 \sim C(x)^3 &\geq& A(x)B(x)C(x) \sim 
|\hat A||\hat B||\hat C|\nonumber \\
&\geq& |\hat A + \hat B + \hat C| \gg {x \over \log^\kappa x},\nonumber
\end{eqnarray}
Thus, $|\hat C| > \kappa$ for $x$ sufficiently large; and, 
$|\hat A + \hat B| > \kappa$, by (\ref{hat_ABC_inequalities}).  
It follows then that $|F|, |G| > \kappa$ for $x$ sufficiently large.

We also have that $F+G \subseteq (\sqrt{x},\infty)$,
since $\hat A, \hat B, \hat C$ satisfy the conclusion to Proposition \ref{extract_B_C}.
Thus, $F$ and $G$ satisfy the hypotheses, and therefore the conclusion, 
of Proposition \ref{elsholtz_prop} with $\delta = \kappa$.  Thus,

\begin{equation}\label{BC_bounds}
{\sqrt{x} \over \log^{6+\kappa}x} \ll |\hat A + \hat B|, |\hat C| \ll 
{\sqrt{x}\log^6x}.
\end{equation}

Between the sets $\hat A$ and $\hat B$, let $S$ be the one with the smaller
number of elements, and let $L$ be the set with the larger number of elements.
We now distinguish two cases:  Case 1 is where $|S| > \kappa$, and Case 2
is where $|S| \leq \kappa$.  

To prove (the contrapositive of) the Main Theorem 
in Case 1, we consider
the two sets $L + \hat C$ and $S$, and let $F$ be the set 
with the smaller number
of elements, and $G$ be the one with the larger number of elements.
(Note:  The sets $F$ and $G$ have now changed from how we defined them
before.)  These sets $F$ and $G$ satisfy the hypotheses of Proposition 
\ref{elsholtz_prop} with $\delta = \kappa$, since 
$$
|S| > \kappa\ {\rm and\ } |L+\hat C| \geq |L| \geq |S| > \kappa\ \Rightarrow\ 
|G| \geq |F| > \kappa,
$$
and since $F$ and $G$ satisfy the other hypotheses of the Proposition.
As in (\ref{BC_bounds}), we deduce from this that 
$$
{\sqrt{x} \over \log^{6+\kappa}x} \ll |L+\hat C|, |S| \ll {\sqrt{x}\log^6x}.
$$ 
From this and (\ref{BC_bounds}) we deduce
$$
A(x), B(x), C(x)\ \gg\ {\sqrt{x} \over \log^{6+\kappa}x}.
$$
Now, since $A,B,C$ is a regular triple, this bound and Lemma
\ref{regular_lemma} give
$$
(A + B + C)(x)\ >\ x^{3/2 - o(1)},
$$
which is absurd.

We now consider Case 2, which is where $|S| \leq \kappa$.
For this case we will have from (\ref{BC_bounds}) that
\begin{equation}\label{LC_inequality}
{\sqrt{x} \over \log^{6+\kappa}x} \ll |L|, |\hat C| \ll \sqrt{x}\log^6 x.
\end{equation}

We need the following Proposition to find subsets of 
$L$ and $\hat C$ with usable properties.

\begin{proposition}\label{regularity_prop}  
There exists a constant $D > 0$ such that if $x$ is sufficiently large,
and if $|S| \leq \kappa$ (Case 2), 
then there exist subsets $L^* \subset L$ and $C^* \subset \hat C$ with
\begin{equation}\label{LC_star}
|L^*|\ >\ {|L| \over \log^D x},\ {\rm and\ } |C^*|\ >\ {|\hat C| 
\over \log^D x},
\end{equation}
such that
\begin{equation}\label{product_sum_inequality}
|L^*|\ |C^*|\ \leq\ 2|L^* + C^*|.
\end{equation}
\end{proposition}

Let $s_1,s_2 \in S$, with $s_2 > s_1$, be any two integers, set 
$k=s_2 - s_1$, and let
$$
L^{\#}\ =\ L^* + s_1\ =\ \{ \ell+s_1\ :\ \ell \in L^*\}.
$$
Then, 
$$
L^{\#} + C^*\ =\ \{\ell+c+s_1\ :\ (\ell,c) \in L^* \times C^*\}
$$
consists entirely of primes, and so does
$$
L^{\#} + C^* + k.
$$
We will need the following Proposition and its Corollary to unlock the 
structure of the set $L^{\#} + C^*$: 

\begin{proposition}\label{few_solutions_prop} 
Let
$$
Q = \max( |L^{\#}|, |C^*|)\ \gg\ {x^{1/2} \over \log^{\kappa+D+6}x},
$$
by (\ref{LC_star}) and (\ref{LC_inequality}).
Then, for any integer $j \geq 1$ we will have
\begin{eqnarray}
&& \sum_{p \leq Q} (\log p)\ \#\{(\ell,c) \in L^{\#} \times C^*\ :\ 
\ell+c+jk \equiv 0 \pmod{p}\} \nonumber \\
&& \hskip1in \ =\ O(j\ |L^{\#}|\ |C^*|\ \log\log x). \nonumber
\end{eqnarray}
\end{proposition}
\bigskip

\begin{corollary}\label{few_solutions_corollary}  
There exists a constant $c >0$ such that all but at most 
$|L^{\#}+C^*|/2$ of the elements $n \in L^{\#}+C^*$ have
\begin{equation}\label{njk_bounds}
\sum_{j=1}^r \sum_{p \leq Q \atop {p | n+jk \atop p\ {\rm prime}}} \log p\ <\ 
cr^2 \log\log x.
\end{equation}
\end{corollary}
\bigskip

One more lemma will establish the Main Theorem:
\bigskip

\begin{lemma}\label{divisor_lemma}  For $x$ sufficiently large, 
there are at most $x\log^{-r/2} x$ integers $n \leq x$ 
which satisfy (\ref{njk_bounds}).
\end{lemma}
\bigskip

We have from Proposition \ref{extract_B_C},  
Corollary \ref{few_solutions_corollary}, Lemma 
\ref{divisor_lemma}, and Proposition \ref{regularity_prop} that for 
$r=4D+2\kappa+2$ and $x$ sufficiently large,
\begin{eqnarray}
&&(A + B + C)(x)\nonumber \\
&&\hskip0.5in  \leq\ 2\ |S\ +\ L\ +\ \hat C|\ \  
{\rm [ Prop.\ \ref{extract_B_C}]} \nonumber \\
&&\hskip0.5in  \leq\ 2\kappa\ |L|\ |\hat C| \nonumber \\
&&\hskip0.5in  <\    2\kappa\ |L^*|\ |C^*|\ \log^{2D}x\ \ 
{\rm [Prop.\ \ref{regularity_prop}]} \nonumber \\
&&\hskip0.5in  \leq\ 4 \kappa\ |L^* + C^*| \log^{2D}x \nonumber \\
&&\hskip0.5in  =\ 4\kappa\ |L^\# + C^*| \log^{2D}x \nonumber \\
&&\hskip0.5in  \leq\ 8 \kappa\ (\log^{2D}x)\ \#\{ n \in L^\# + C^*\ :\ 
n\ {\rm satisfies\ } (\ref{njk_bounds}) \}\ \ 
{\rm [Cor.\ \ref{few_solutions_corollary}]} \nonumber \\
&&\hskip0.5in  \leq\ 8 \kappa\ (\log^{2D}x)\ \#\{ n \leq 3x\ :\ n\ 
{\rm satisfies\ } (\ref{njk_bounds}) \} \nonumber \\
&&\hskip0.5in  \leq\ 8 \kappa\ (\log^{2D}x)\ {3x \over \log^{2D+\kappa+1} 
(3x)} \nonumber \\
&&\hskip0.5in  \ll_\kappa\ {x \over \log^{\kappa+1/2}x}\ \ 
{\rm [Lemma\ \ref{divisor_lemma}]}. \nonumber
\end{eqnarray}
which contradicts (\ref{main_theorem_hypothesis}), 
and so the theorem is proved.
\bigskip

\section{Proof of Theorem \ref{unconditionallowerbound}.}

The proof is based on a twofold application of a counting argument due
to Erd\H{o}s, Stewart and Tijdeman \cite{ErdosandStewartandTijdeman},
compare also Lemma 6 in Pomerance, S\'ark\"ozy and Stewart 
\cite{PomeranceandSarkozyandStewart}.
\begin{lemma}
Let $\tau$ be a positive integer.
Let $x> x_{\tau} $ be a sufficiently large positive integer and let 
$T$ be a non-empty subset of $\{1, \ldots , x\}$.  
Then there exists
a set $S\subset T$ and a set of non-negative integers $A$ such that
$A+S \subset T$, and 
$$
|S| \geq {{|T| \choose \tau} \over \binom{x-1}{\tau -1}},\ |A|=\tau.
$$
\end{lemma}
Since we want to prescribe the number of elements
of two sets $A$ and $B$ we apply this lemma once again to the set 
$S$. 
This gives
\begin{corollary}
Let $\kappa_1, \kappa_2$ denote positive integers.
Let $x > x_{\kappa_1,\kappa_2}$ be a sufficiently large positive integer 
and let $T$ be a non-empty subset of $\{1, \ldots , x\}$. 
Let 
$$
R= {\binom{|T|}{\kappa_1} \over \binom{x-1}{\kappa_1-1}}.
$$
Then there exists a subset $C \subset T$ and sets of non-negative 
integers $A,B$ such that
$$ 
A+B+C \subset T,\ |C| \geq {\binom{R}{\kappa_2} \over \binom{x-1}{\kappa_2-1}},
\ |A|=\kappa_1,\quad |B|=\kappa_2. 
$$
\end{corollary}
It is obvious that one could iterate this argument.
We resist doing this since we concentrate on ternary problems.

Now let $T$ denote the set of primes in $[1,x]$. Recall that
by the prime number theorem with error term 
(see \cite{Landau}, \textsection 54)
$$
|T|= \frac{x}{\log x}+ \frac{x}{(\log x)^2} + 
O\left( \frac{x}{(\log x)^3}\right).
$$
For large $x$ we have that $|T|- \kappa_1> \frac{x}{\log x}$.
Hence it follows (as in the proof of theorem 6 in
\cite{PomeranceandSarkozyandStewart})
that 
\begin{eqnarray}
R&\geq & \frac{\frac{1}{\kappa_1!}
 \left( \frac{x}{\log x} + \frac{x}{2(\log x)^2}\right)^{\kappa_1} 
}{\frac{1}{(\kappa_1-1)!} x^{\kappa_1-1}}
\geq \frac{x}{\kappa_1 (\log x)^{\kappa_1}} + 
\frac{x}{2 (\log x)^{\kappa_1 +1}}
\end{eqnarray}
For the second application of the argument we observe that for large $x$ we
have
$R-\kappa_2> \frac{x}{\kappa_1 (\log x)^{\kappa_1}}$.
The argument then gives:
\begin{eqnarray}
|C|&\geq & \frac{ 
\frac{1}{\kappa_2!} 
\left(\frac{x}{\kappa_1 (\log x)^{\kappa_1}}\right)^{\kappa_2}}{
\frac{1}{(\kappa_2 -1)!x^{\kappa_2-1}}}
= \frac{x}{\kappa_2 \kappa_1^{\kappa_2} (\log x)^{\kappa_1 \kappa_2}}
\end{eqnarray}
Our theorem now follows since $|C| \leq |A+B+C|$.

\section{Statements and Proofs of Some Technical Lemmas.} 
\label{technical_section}

We will need the following three sieve lemmas, and their various 
corollaries:  the Large Sieve of Montgomery (see \cite{Montgomery}), 
Brun's Upper Bound Sieve (see \cite{Halberstam}), and
a variant of Gallagher's Sieve (see \cite{Gallagher}):

\begin{lemma}{\bf (Montgomery's Sieve)} \label{large_sieve}
Given a set of integers $J \subseteq \{1,2,...,x\}$,
and for each prime $p \leq x$, let $\omega(p)$ be the number of progressions
modulo $p$ which $J$ fails to occupy.  Then,
$$
|J|\ \leq\ {x+Q^2 \over \sum_{q \leq Q} \mu^2(q) \prod_{p | q} 
{\omega(p) \over p-\omega(p)} }.
$$
\end{lemma}
\bigskip

One has the following corollary, which essentially appears 
in Vaughan's paper \cite{Vaughan}.

\begin{corollary} \label{large_sieve_corollary} 
For $J$ and $\omega(p)$ as above, and $T \leq \sqrt{x}$, we have
$$
|J|\ \leq\ 
{2x \over \left ({1 \over m} \sum_{p \leq T} {\omega(p) \over p} \right )^m},
$$
where $m = \lfloor (\log x)/(2\log T) \rfloor$.
\end{corollary}
\noindent (Note:  In Vaughan's paper he proves this result with the factor $4$
on the left hand side, instead of the factor $2$.  The reason is that he
used an earlier, weaker form of the Large Sieve.)


\begin{lemma}\label{bruns_sieve} {\bf (Brun's Sieve)}
Suppose that $J \subseteq \{1,2,...,x\}$ is the largest such set of integers
which fails to occupy $\omega(p) \leq B$ progressions modulo $p$, for
each prime $p \leq z$.  Then,
$$
|J|\ \ll_B\ x\ \prod_{p \leq z} \left (1 - {\omega(p) \over p} \right ).
$$
\end{lemma}

\begin{lemma}{\bf (Gallagher's Sieve)}\label{gallagher_sieve1} 
Suppose that $J \subseteq \{1,2,...,x\}$, and
$|J| > U$.  Then,
$$
|J|^2 (\log x + O(1))\ >\ \sum_{p \leq U\atop p\ {\rm prime}} 
(\log p) \sum_{c=0}^{p-1} |J(c,p)|^2.
$$
\end{lemma}
\bigskip

A corollary of this sieve which we will need is the following:

\begin{corollary}\label{gallagher_sieve2}
Suppose $J$ is as in Lemma \ref{gallagher_sieve1}, and let $h(p)$
denote the number of residue classes modulo $p$ occupied by $J$, for each
$p \leq U < |J|$.  Then,
\begin{eqnarray}
\log x + O(1)\ &>&\ {1 \over |J|^2} \sum_{p \leq U} (\log p) \sum_{c=0}^{p-1}
|J(c,p)|^2 \nonumber \\
&\geq&\ \sum_{p \leq U} {\log p \over h(p)}.\nonumber
\end{eqnarray}
\end{corollary}
\bigskip

We will also need the following inequality of Cauchy and Davenport (see \cite{Nathanson}):

\begin{lemma}\ {\bf (Cauchy-Davenport Inequality)}\label{cauchy_davenport_lemma}
For sets $G$ and $H$,
let $h_1, h_2$ and $h_3$ denote the number of residue classes modulo $p$
occupied by $G, H$ and $G+H$, respectively.  Then,
$$
h_3\ \geq\ \min(h_1+h_2-1, p).
$$
\end{lemma}
\bigskip

Finally, we will also need the following simple consequence of the Cauchy-Schwarz
inequality:

\begin{lemma}\label{square_lower_bound}  Suppose that $J$ is a set of 
integers which occupies at most $k$ progressions modulo $m$.  Then,
$$
\sum_{a=0}^{m-1} |J(a,m)|^2\ \geq\ {|J|^2 \over k}.
$$
\end{lemma}

To prove this lemma, let $\delta(a)$ be $1$ if $a$ is in
one of the progressions occupied by $J$ (there are at most $k$ such
progressions), and let it be $0$ otherwise.  
Then, the Lemma follows quickly from the Cauchy-Schwarz inequality:
\begin{eqnarray}
k \sum_{a=0}^{m-1} |J(a,m)|^2 &=&
\left ( \sum_{a=0}^{m-1} |J(a,m)|^2 \right )\left (
\sum_{a=0}^{m-1} \delta(a)^2 \right )\nonumber \\
&\geq&\ \left ( \sum_{a=0}^{m-1} |J(a,m)|\delta(a)\right )^2\ =\ 
|J|^2.\nonumber
\end{eqnarray}

We now prove those of the above lemmas which cannot be found in the
literature, as well as Lemmas \ref{regular_lemma} and \ref{divisor_lemma}.
\bigskip

\noindent {\bf Proof of Lemma \ref{regular_lemma}.}  
To prove this lemma we let $\epsilon = 1/2$, and let $D > 0$ 
and $S$ be as in the definition of regular sets.  Further,
let 
$$
T = \{ n \in A_1 + \cdots + A_k\ :\ n \leq x,\ r(n;A_1,...,A_k) \leq 
\log^Dx\}.
$$
Then, 
\begin{eqnarray}
A_1(x) \cdots A_k(x) &=& 
\sum_{n \in T} r(n;A_1,...,A_k)\nonumber \\
&&\ \ \ \ \ \  + \sum_{n \in S}r(n;A_1,...,A_k) \nonumber \\
&<& (\log^D x)T(x) + {1 \over 2} A_1(x) \cdots A_k(x)
\nonumber \\
&\leq& (\log^D x)\ (A_1+\cdots+A_k)(x)\nonumber \\
&&\ \ \ \ \ \  + {1 \over 2} A_1(x) \cdots A_k(x). \nonumber
\end{eqnarray}
Rearranging terms gives
$$
A_1(x) \cdots A_k(x)\ <\ 2(\log^D x)\ (A_1+ \cdots + A_k)(x);
$$
and so, the conclusion of the Lemma holds with $E = D+1$.
\bigskip

\noindent {\bf Proof of Lemma \ref{divisor_lemma}.}
We note that if $n$ satisfies (\ref{njk_bounds}) then the largest prime 
divisor of each of the numbers $n+k,n+2k,...,n+rk$ is $< \log^{cr^2} x$.
Thus, for each prime $p \in [\log^{cr^2} x, x)$, we must have that 
$$
n \not \equiv -k, -2k, ..., -rk\ \pmod{p}.
$$
Thus, the number of progressions which $n$ can lie in modulo $p$, for
each such $p$, is 
$h(p) < p-r$.   From Brun's Sieve, we get that the number of integers
$n$ satisfying (\ref{njk_bounds}) is
\begin{eqnarray}
\ll_r\ x\ \prod_{\log^{cr^2} x < p < Q\atop p\ {\rm prime}} 
\left ( 1 - {r \over p} \right )
&<&\ x\ \exp \left ( - r\ \sum_{\log^{cr^2} x < p < Q \atop p\ {\rm prime}}
{1 \over p} \right ) \nonumber \\
&=&\ x\ \exp \biggl ( -r (\log\log x\ -\ O(\log\log\log x)) \biggr ) 
\nonumber \\ 
&=&\ o \left ( {x \over \log^{r/2} x} \right ),
\end{eqnarray}
which proves the Lemma.
\bigskip

\noindent {\bf Proof of Lemma \ref{gallagher_sieve1}.}
\bigskip
We have that for any pair of integers $j_1,j_2 \in J$, $|j_1-j_2| < x$, and
so 
$$
\sum_{p | j_1-j_2\atop p\ {\rm prime}} \log p\ \leq \log |j_1-j_2|\ <\ \log x.
$$
Summing over all pairs $j_1,j_2$ (of which there are at most $|J|^2$), we get
\begin{eqnarray}
|J|^2 \log x\ &>&\ 
\sum_{j_1,j_2 \in J\atop j_1\neq j_2} \sum_{p | j_1-j_2 \atop p\ {\rm prime}}
\log p\nonumber \\
&>&\ \sum_{p \leq U} (\log p) \sum_{c=0}^{p-1} \#\{j_1,j_2\in J\ :\ j_1 \neq j_2,
\nonumber \\
&&\hskip1.5in\ j_1 \equiv j_2 \equiv c \pmod{p}\}\nonumber \\
&=&\ \sum_{p \leq U} \left (\ (\log p)\left 
(\sum_{c=0}^{p-1} |J(c,p)|^2 \right )\ -\ |J|\log p \right ).\nonumber
\end{eqnarray}
Using the fact that
$$
\sum_{p \leq U} \log p\ =\ O(U)\ =\ O(|J|),
$$
and rearranging terms in the above string of inequalities, we get
$$
|J|^2 (\log x + O(1))\ >\ \sum_{p \leq U} (\log p) \sum_{c=0}^{p-1}
|J(c,p)|^2,
$$
as claimed.
\bigskip

\noindent {\bf Proof of Corollary \ref{gallagher_sieve2}.}
\bigskip

Since $J$ occupies $h(p)$ progressions modulo $p$, we have
from Lemma \ref{square_lower_bound} that
$$
{1 \over |J|^2} \sum_{c=0}^{p-1} |J(c,p)|^2\ >\ {1 \over h(p)}.
$$

Putting this into Lemma \ref{gallagher_sieve1}, we get
\begin{eqnarray}
\log x +O(1))\ &>&\ {1 \over |J|^2} \sum_{p \leq U \atop p\ {\rm prime}}
(\log p) \sum_{c=0}^{p-1} |J(c,p)|^2\ >\ \sum_{p \leq U} {\log p \over h(p)},
\nonumber
\end{eqnarray}
as claimed.
\bigskip

\section{Proof of Proposition \ref{elsholtz_prop}.}

Let $\tau = \lfloor \delta \rfloor+1$.  Then, we have
$|G| \geq |F| \geq \tau$.  

The proof involves four iterations:  In the first iteration we will 
show that $|G| \ll x(\log x)^{-\tau+o(1)}$,
and thus $|F| \gg (\log x)^{\tau-\delta-o(1)}$; in the second 
iteration, we will show that 
$|F| \gg \exp( (\log x)^{\tau-\delta})$; in the third iteration, we 
will show that $|F| > x^{1/3}$, for $x$ sufficiently large; 
and, in the final iteration, we will show that 
$$
\frac{x^{1/2}}{\log^{\delta+6}} \ll |F| \leq |G| \ll x^{1/2}\log^6x.
$$

We note that our Proposition can be proved using three iterations
(instead of four), as was done in \cite{Elsholtz2}; also, no attempt
was made to optimize the powers of the logarithms appearing in the result. 
 

Throughout the proof we let $h_1(p)$ and $h_2(p)$ denote the number 
of residue classes occupied by $F$ and $G$, respectively.  Since no
element of $F+G$ can be divisible by a prime $\leq \sqrt{x}$, we deduce
that $F+G$ occupies at most $p-1$ residue classes modulo $p$ for each
such prime.  So, from Lemma \ref{cauchy_davenport_lemma}, we deduce
\begin{eqnarray}\label{cauchy_davenport_conclusion}
h_1(p) + h_2(p) \leq p.
\end{eqnarray}
We let $\omega(p) = p-h_2(p)$ be the number of progressions which $G$ 
fails to occupy; and so, (\ref{cauchy_davenport_conclusion}) implies that 
$\omega(p) \geq h_1(p)$.  

For the first iteration, let $f_1,...,f_\tau$ be any $\tau$ elements of
$F$, and $Z$ be the set of primes $\leq \sqrt{x}$ with the property 
that $f_1,...,f_\tau$ all occupy different residue classes modulo $p$.
Let $P$ be the set of primes $\leq \sqrt{x}$, and set
$$
f(Z) = \sum_{p \in Z} {1 \over p} = \log\log x + O(1) - \sum_{p \in
P \setminus Z} {1 \over p}.
$$
To estimate this last sum, we first define 
$$
s(n) = \sum_{p | n\atop p\ {\rm prime}} {1 \over p}.
$$
Then, $s(n) \ll \log\log\log n$, and this upper bound is attained
when $n$ is the product of the primes $\leq \log n$.  
Now, if $p \in P \setminus Z$, then $p | \Delta$, where
$$
\Delta = \prod_{1 \leq i < j \leq \tau} |f_j - f_i| \ll x^{\tau^2/2};
$$
and so,
$$
\sum_{p \in P \setminus Z} {1 \over p} \leq \sum_{p | \Delta} {1 \over p}
\ll_\tau \log\log\log x.
$$
Thus,
$$
f(Z) = \log\log x - O_\tau(\log\log\log x).
$$
Letting $\Pi(x)$ be the product of the primes $\leq \sqrt{x}$,
we deduce from Lemma \ref{bruns_sieve} that 
\begin{eqnarray}
|G| &\leq& \#\{n \leq x\ :\ ((n+f_1)(n+f_2)\cdots (n+f_\tau),\Pi(x))=1\}\nonumber \\
&\ll_\tau& x \prod_{p \in Z} \left ( 1 - {\tau \over p} \right )\nonumber \\
&\ll& x \exp \left ( -\tau \sum_{p \in Z} {1 \over p} \right )\nonumber \\
&<& {x \over \log^{\tau -o(1)}x}.\nonumber
\end{eqnarray}
Thus, since $x(\log x)^{-\delta} \ll |F||G|$, we deduce
$|F| \gg (\log x)^{\tau - \delta -o(1)}$, as claimed.

For the second iteration, let $f_1,...,f_t \in F$, where
$t = \log^{\tau-\delta - o(1)}x$, and, as before, let
$Z'$ be the set of primes $\leq \sqrt{x}$ where all the
$f_i$'s fall into distinct residue classes modulo $p$.  
Then, as before, let 
$$
\Delta' = \prod_{1 \leq i < j \leq t} |f_i - f_j| \ll x^{t^2/2}. 
$$
Then, $f_1,...,f_t$ are not distinct modulo $p$ implies $p | \Delta'$.
As before, we have
$$
\sum_{p \in P \setminus Z'} {1 \over p} \leq \sum_{p | \Delta'} {1 \over p}
\ll \log\log\log x.
$$
Thus, if we let $T = \exp( \log^{1 - \tau/2 + \delta/2}x)$, then
\begin{eqnarray}
\sum_{p \leq T \atop p \in Z'} {1 \over p}\ &\geq&\ \sum_{p \leq T\atop p\ 
{\rm prime}} {1 \over p}\ -\ \sum_{p \in P \setminus Z'} {1 \over p}
\nonumber \\ 
&=& \log\log T\ -\ O(\log\log\log x). \nonumber 
\end{eqnarray}
Now, applying Corollary \ref{large_sieve_corollary} with 
$$
m = \lfloor (\log^{\tau/2 - \delta/2}x)/2\rfloor\ {\rm and\ } 
\omega(p) = (\log x)^{\tau-\delta-o(1)}\ {\rm for\ all\ }p \in Z' \cap [2,T],
$$ 
we get
\begin{eqnarray}
|G| &\leq& \#\{n \leq x\ :\ ((n+f_1)(n+f_2)\cdots (n+f_t),\Pi(x))=1\}\nonumber \\
&\leq&  {x \over \left ( (\log x)^{\tau/2 - \delta/2 - o(1)}\sum_{p\leq T
\atop p \in Z'} {1 \over p} \right )^m} \nonumber \\
&\ll& {x \over \exp(2m)} \ll {x \over \exp(\log^{\tau/2 - \delta/2}x)}. \nonumber
\end{eqnarray}
Thus, since $x/\log^\delta x < |F||G|$, we conclude
that 
$$
|F| > \exp \left ((1-o(1))\sqrt{\log^{\tau-\delta}x} \right ).
$$ 

For the third iteration, let $T' = \exp(\sqrt{\log^{\tau-\delta} x}/2)$
and $m' = \lfloor \sqrt{\log^{2-\tau+\delta}}\rfloor$.
Then, from Corollary \ref{gallagher_sieve2}, we have
$$
\log x + O(1) > \sum_{T'/2 \leq p \leq T'} {\log p \over h_1(p)}.
$$
(Note: we use the corollary with $J=F$, and we have from iteration two
that $|F| > T'$ for $x$ sufficiently large).
So, for almost all primes $p \in [T'/2, T']$ we have that 
$\omega(p) \geq h_1(p) > p/\log^2x$; 
and so, 
$$
\sum_{T'/2 \leq p \leq T} {\omega(p) \over p} \gg {T' \over \log^3x}.
$$
Using Corollary \ref{large_sieve_corollary} with $T=T'$, we deduce 
$$
|G| \ll {x \over \left ({1 \over m'} \sum_{T'/2 \leq p \leq T'} {\omega(p) \over p} \right )^
{m'}} < {x \over x^{1/2 - o(1)}} = x^{1/2+o(1)}.
$$
Thus, since $|F||G| \gg x\log^{-\delta}x$, we deduce $|F| > x^{1/3}$
for $x$ sufficiently large.

For the last iteration, we have by Corollary \ref{gallagher_sieve2}
that
$$
\log x + O(1) > \sum_{x^{1/4}/2 \leq p \leq x^{1/4}} {\log p \over h_1(p)};
$$
and it follows that almost all primes in $[x^{1/4}/2,x^{1/4}]$
have $\omega(p) \geq h_1(p) > p/\log^2 x$.  Thus,
$$
\sum_{x^{1/4}/2 \leq p \leq x^{1/4}} {\omega(p) \over p} \gg {x^{1/4}
\over \log^3x}.
$$ 
By Corollary \ref{large_sieve_corollary} we have
$$
|G| \ll\ {x \over \left ( \sum_{x^{1/4}/2 \leq p \leq x^{1/4}} 
{\omega(p) \over p} \right )^2} \ll {x\log^6x \over x^{1/2}} 
= \sqrt{x}\log^6x;
$$
and so, since $x/\log^\delta x \ll |F||G|$, we deduce
$$
{\sqrt{x} \over \log^{6+\delta}x} \ll |F| \leq |G| \ll {\sqrt{x}\log^6x},
$$ 
and the Proposition is proved.
\bigskip

\section{Proof of Proposition \ref{extract_B_C}.}

Since $A,B,C$ is a regular triple of sets, we have from Lemma 
\ref{regular_lemma} that for some $E > 0$,
\begin{eqnarray}
{A(\sqrt{x}) \over A(x)} \cdot {B(\sqrt{x}) \over B(x)} 
\cdot {C(\sqrt{x}) \over C(x)} &\leq& (\log^E x)\ {(A+B+C)(\sqrt{x})
\over (A+B+C)(x)} \nonumber \\
&\ll& {\sqrt{x}\log^E x \over x/\log^\kappa x} = {\log^{E+\kappa} x 
\over \sqrt{x}}.
\end{eqnarray}
Thus, for $x$ sufficiently large, one of the following inequalities must 
hold:
\begin{eqnarray}\label{either_or_bounds}
A(\sqrt{x}) < {A(x) \over x^{1/5}}\ \ {\rm or\ \ }
B(\sqrt{x}) < {B(x) \over x^{1/5}}\ \ {\rm or\ \ }
C(\sqrt{x}) < {C(x) \over x^{1/5}}.
\end{eqnarray}
Suppose that the inequality holds for $A(x)$ and $A(\sqrt{x})$.  Then, 
letting 
$$
\hat A = A \cap (\sqrt{x},x],\ \hat B = B \cap [1,x],\  
{\rm and\ }\hat C = C \cap [1,x]
$$
gives
$$
|\hat A| \sim A(x),\ |\hat B| \sim B(x),\ |\hat C| \sim C(x),\ 
{\rm and\ } |\hat A + \hat B + \hat C| \subset (\sqrt{x},x].
$$
Also, since $A,B,C$ is a regular triple, we get
\begin{eqnarray}
0 &\leq& (A+B+C)(x) - |\hat A + \hat B + \hat C| \nonumber \\
&\leq& 
\#\{n = a+b+c\ :\ a\in A, b\in B, c \in C,\ a \leq \sqrt{x}\} \nonumber \\
&\leq& A(\sqrt{x})\ B(x)\ C(x) \nonumber \\
&\leq& x^{-1/5} A(x)\ B(x)\ C(x) \nonumber \\
&\leq& x^{-1/5}(\log^E x)\ (A+B+C)(x) \nonumber
\end{eqnarray}
Thus,
$$
(A+B+C)(x)\ \sim\ |\hat A + \hat B + \hat C|,
$$
as claimed.  We get the same conclusions for the remaining cases of
(\ref{either_or_bounds}).
\bigskip

\section{Proof of Proposition \ref{regularity_prop}.}

Since $A,B,C$ is a regular triple, one can easily deduce that for
$\epsilon = 1/12$ and $x$ sufficiently large, there exists $E > 0$
such that if $|S| \leq \kappa$, then
\begin{equation}\label{regular_LC_bounds}
\sum_{n \in L \times \hat C\atop r(n; L, \hat C) > \log^E x} r(n; L, \hat C)
\ <\ \epsilon |L \times \hat C|.
\end{equation}
For the remainder of the proof of this Proposition, we will assume 
that $E$ is such that this inequality is satisfied.

The proof now proceeds using a probabilistic argument:  Let 
$L'$ and $C'$ be random subsets of $L$ and $\hat C$, respectively, where
$$
{\rm Prob}(\ell \in L'\ |\ \ell \in L)\ =\ 
{\rm Prob}(c \in C'\ |\ c \in \hat C)\ =\ {1 \over \log^{2E}x},
$$
where all these probabilities are independent.
Clearly, $|L'|$ and $|C'|$ each have a binomial distribution, 
which implies that the following occurs with probability $1- o(1)$:

\begin{equation}\label{central_inequality}
{E(|L' \times C'|) \over 2}\ <\ |L' \times C'|\ <\ 
2 E(|L' \times C'|),
\end{equation}
where $E(|L' \times C'|)$ is  the usual expectation given by
$$
E(|L' \times C'|)\ =\ \sum_{(\ell,c) \in L \times \hat C} {\rm Prob}
((\ell,c) \in L' \times C')\ =\ {|L \times \hat C| \over \log^{4E} x}.
$$ 

In the course of our proof, we will show that the event

\begin{equation}\label{random_inequality}
(1-6\epsilon)|L' \times C'|\ <\ |L' + C'|\ \ {\rm and }\ \  
(\ref{central_inequality})\ {\rm occurs}
\end{equation}
has positive probability, which will imply that there exists
subsets $L^* \subset L$ and $C^* \subset \hat C$ satisfying
these same inequalities.  If we can do this, then 
(\ref{product_sum_inequality}) will hold (since $\epsilon = 1/12$),
and (\ref{LC_star}) will hold for $D = E+1$ and 
$x$ sufficiently large.

Thus, the Proposition will follow if we can show that 
(\ref{random_inequality}) has positive probability.  
We note that it suffices to prove that
\begin{equation}\label{target_inequality}
{\rm Prob}\left (\ |L' \times C'| - |L' + C'|\ <\ 3\epsilon E(L' \times C')
\ \right )\ >\ {1 \over 2},
\end{equation}
since (\ref{central_inequality}) holds with probability $1-o(1)$.

We now proceed to show that (\ref{target_inequality}) holds:
Suppose that $n \in L + \hat C$ has exactly $k$ solutions to 
$$
n\ =\ \ell_1 + c_1, ..., \ell_k + c_k,\ {\rm each\ }(\ell_i,c_i) 
\in L \times \hat C.
$$
Then, since the $\ell_i$'s are distinct, and the $c_i$'s distinct,
we have that all subsets of the following probabilities are independent
$$
{\rm Prob}((\ell_1,c_1) \in L' \times C'),..., {\rm Prob}((\ell_k,c_k) \in L' 
\times C')\ =\ {1 \over \log^{4E}x}.
$$
It follows then that if we let $r'(n)$ be the random variable
$$
r'(n)\ =\ \{(\ell,c) \in L' \times C'\ :\ n = \ell+c\},
$$
then ${\rm Prob}(r'(n)=d)$ has a binomial distribution, given by
$$
{\rm Prob}(r'(n) = d) = {k \choose d}\left ( 1 - {1 \over \log^{4E}x}
\right )^{k-d} {1 \over \log^{4dE} x} < {k^d \over d! \log^{4dE}x};
$$
and, we have the following easily checked expectation formula
$$
E(r'(n))\ =\ {r(n;L,\hat C) \over \log^{4E}x},
$$
where $r(n;L,\hat C)$ is as defined in the Introduction.

For bookkeeping purposes, define
\begin{eqnarray}
N &=& \{ n \in L + \hat C\ :\ n\ {\rm has\ at\ most\ }\log^Ex\ {\rm solutions\ to}
\nonumber \\
&&\hskip1in n = \ell+c,\ (\ell,c) \in L \times \hat C \};{\rm and\ } \nonumber \\
\overline{N} &=& (L + \hat C)\ \setminus\ N. \nonumber 
\end{eqnarray}
and define the random variable
$$
\delta(n) = \begin{cases} 0, &{\rm if}\ n \not \in L' + C'; \\
                          1, &{\rm if}\ n \in L'+C'.\end{cases}
$$
Then, from (\ref{regular_LC_bounds}) and the above probability and
expectation estimates, we have:
\begin{eqnarray}
&& E( |L' \times C'| - |L' + C'| )\nonumber \\
&&\hskip0.5in =\ \sum_{n \in L + \hat C}
E(r'(n) - \delta(n)) \nonumber \\
&&\hskip0.5in =\ \sum_{n \in N} E(r'(n) - \delta(n)) + \sum_{n \in \overline{N}}
E(r'(n) - \delta(n)) \nonumber \\
&&\hskip0.5in \leq\ \sum_{n \in N} \sum_{d \geq 2} (d-1)\ {\rm Prob}(r'(n)=d)
\ +\ \sum_{n \in \overline{N}} E(r'(n)) \nonumber \\
&&\hskip0.5in \leq\ \sum_{n \in N} \sum_{d \geq 2} {1 \over (d-1)!\log^{3dE}x} 
\ +\ {1 \over \log^{4E}x}\sum_{n \in \overline{N}}
r(n; L, \hat C) \nonumber \\
&&\hskip0.5in \leq\ {2|L \times \hat C| \over \log^{6E}x}\ +\ {\epsilon\ |L \times \hat C|
\over \log^{4E}x}\ =\ E(|L' \times C'|) 
\left ( \epsilon + {2 \over \log^{2E}x} \right ). \nonumber
\end{eqnarray}
\bigskip

\noindent {\bf Markov's Inequality.}  If $X$ is a non-negative random variable,
then 
$$
{\rm Prob}(X \geq a)\ \leq\ {E(X) \over a}
$$
\bigskip

From this inequality with $X = |L'\times C'| - |L'+C'|$, together
with our above expectation estimates, we deduce
$$
{\rm Prob} \bigl ( |L' \times C'|\ -\ |L' + C'|\ \geq\  3\epsilon
E(|L' \times C'|) \bigr )\ <\ {1 \over 3} - {3 \over \log^{4E}x}.
$$
Therefore, (\ref{target_inequality}) holds for $x$ sufficiently large.
\bigskip

\section{Proof of Proposition \ref{few_solutions_prop}, Corollary 
\ref{few_solutions_corollary}, and Lemma \ref{cauchy_davenport_consequence}.}

For a given set of integers $J$, let $J_p$ denote the set of 
residue classes modulo $p$ occupied by $J$, and let $\overline{J_p}$
denote those residue classes {\it not} occupied by $J$.  Clearly,
$|\overline{J_p}| = p - |J_p|$. 

For all integers $j$ we have that $(\ell,c) \in L^\# \times C^*$ is
a solution to $\ell+c+jk \equiv 0 \pmod{p}$ if and only if
$$
(\ell,c) \equiv (r,-r-jk) \pmod{p},\ 
{\rm for\ some\ } r\in L^\#_p \setminus (L^\#-jk)_p.
$$

From this and Cauchy's inequality we have 
\begin{eqnarray}
Z &:=& \sum_{p \leq Q} (\log p) \#\{(\ell,c) \in L^{\#} \times C^*\ :\ 
\ell+c+jk \equiv 0 \pmod{p}\} \nonumber \\
&=&\ \sum_{p \leq Q} (\log p) \sum_{r \in L^\#_p \setminus (L^\# - jk)_p}
|L^\#(r,p)|\ |C^*(-r-jk,p)| \nonumber \\
&\leq& Z_1^{1/2} Z_2^{1/2}, \nonumber
\end{eqnarray}
where
\begin{eqnarray}
Z_1 &=& \sum_{p \leq Q} (\log p) \sum_{r \in L^\#_p \setminus (L^\#-jk)_p}
|L^\#(r,p)|^2;\ {\rm and} \nonumber \\
Z_2 &=& \sum_{p \leq Q} (\log p) 
\sum_{r \in L^\#_p \setminus (L^\#-jk)_p} |C^*(-r-jk,p)|^2. \nonumber 
\end{eqnarray}

To bound $Z_1$ and $Z_2$ from above we will require the following three
results:

\begin{lemma}\label{cauchy_davenport_consequence}
We have for $Q = \sqrt{x} \log^{O(1)}x$ that
\begin{eqnarray}
\sum_{p \leq Q} {\log p \over |L^\#_p|} &=& \log x\ +\ O(\log\log x)
\ =\ \sum_{p \leq Q} {\log p \over |C^*_p|};\ {\rm and\ }
\label{LCQ_sum1} \\
\sum_{p \leq Q} {\log p \over p - |L^\#_p|} &=& \log x\ +\ O(\log\log x)
\ =\ \sum_{p \leq Q} {\log p \over p - |C^*_p|}. 
\label{LCQ_sum2}
\end{eqnarray}
\end{lemma}

Compare this with Elsholtz \cite{Elsholtz3}.

\begin{proposition}\label{variance_prop}
Suppose that $J = L^\#$ or $C^*$, and let $K$ be the other set
(if $J = L^\#$, then $K=C^*$, and vice versa).
Also, suppose that for each prime $p \leq Q$ we have a set of residue classes 
$G_p \subseteq J_p$.  Then, we have the following inequality
$$
\sum_{p \leq Q} (\log p) \sum_{r \in G_p} 
|J(r,p)|^2\ <\ |J|^2\left ( \sum_{p \leq Q} {(\log p) |G_p| \over (p-|K_p|)^2}\ 
+\ O(\log\log x) \right ).
$$
\end{proposition}

\begin{lemma}\label{translation_invariant_lemma}
Suppose $J = L^\#$ or $C^*$, and that $K$ is the other of the two sets.
Then, for any integer $j > 0$, 
$$
\sum_{p \leq Q} {(\log p)|J_p\ \setminus\ (J-jk)_p| \over p - |K_p|^2}
\ =\ O(j\ \log\log x).
$$
\end{lemma}

The proofs of these last two results will make use of the following basic
facts about the sets $L^\#$ and $C^*$:  
Since for every $(\ell,c) \in L^\# \times C^*$ we have 
$\ell + c$ and $\ell+c+k$ are primes $> \sqrt{x}$, there can 
be no solutions to $\ell+c \equiv 0 \pmod{p}$ or $\ell+c+k \equiv 0 \pmod{p}$
for any prime $p \leq Q < \sqrt{x}$.  Thus,
\begin{eqnarray}\label{L_subset_C}
&& L^\#_p \cap (-C^*)_p\ =\ \emptyset\ =\ (L^\#+k)_p \cap (-C^*)_p 
\nonumber \\
&& \ \ \ \Longrightarrow\  L^\#_p\ {\rm and\ }\ (L^\#+k)_p\ {\rm are\ 
both\ subsets\ of\ } \overline{(-C^*)_p}. 
\end{eqnarray}
Similarly,
\begin{equation} \label{C_subset_L}
C^*_p\ {\rm and\ }(C^*+k)_p\ {\rm are\ both\ subsets\ of\ }\overline{(-L^\#)_p}.
\end{equation}

Resuming the proof of our Proposition \ref{few_solutions_prop}, 
we have from Proposition
\ref{variance_prop} and Lemma \ref{translation_invariant_lemma} 
with $J = L^\#$ and $G_p = L^\#_p\ \setminus\ (L^\#-jk)_p$ that
$$
Z_1\ <\ |L^\#|^2\left ( \sum_{p \leq Q} {(\log p) |G_p| \over 
(p - |K_p|)^2}\ +\ O(\log\log x) \right )\ =\ O(j\ |L^\#|^2\ \log\log x).
$$ 
Applying these two results with $J = C^*$ and 
$G_p = (- (L^\#\ \setminus\ (L^\#-jk)_p)\ -\ jk)_p$\ (note: 
$|G_p| = |L^\#\ \setminus\ (L^\#-jk)_p|$ ), we likewise get 
$$
Z_2\ =\ O(j\ |C^*|^2\ \log\log x). 
$$
Thus,
$$
Z\ \leq\ Z_1^{1/2} Z_2^{1/2}\ =\ O(j\ |L^\#|\ |C^*|\ \log\log x),
$$ 
which proves the Proposition.

\bigskip    

\noindent {\bf Proof of Corollary \ref{few_solutions_corollary}.}
\bigskip
We have that
\begin{eqnarray}
\sum_{n \in L^\# + C^*} \left (\sum_{j=1}^r \sum_{p \leq Q\atop {p | n+jk
\atop p\ {\rm prime}}} \log p \right ) &\leq& 
\sum_{(\ell,c) \in L^\# \times C^*} 
\sum_{j=1}^r \sum_{p \leq Q \atop { p | \ell+c+jk \atop p\ {\rm prime}}}
\log p \nonumber \\
&=& \sum_{j=1}^r \sum_{p \leq Q} \#\{ (\ell,c) \in L^\# \times C^*\ :\ 
p\ |\ \ell+c+jk \} \nonumber \\
&=& \sum_{j=1}^r O(j\ |L^\#|\ |C^*|\ \log\log x) \nonumber \\
&=& O(r^2\ |L^\# + C^*|\ \log\log x), \nonumber 
\end{eqnarray}
where this last equality following from Proposition \ref{regularity_prop}.
It is now obvious that more than half the elements 
$n \in L^\# + C^*$ satisfy (\ref{njk_bounds}), which proves the Corollary. 
\bigskip

\noindent {\bf Proof of Lemma \ref{cauchy_davenport_consequence}.}  

Since for $p \leq Q$ we have $L^\#+C^*$ contains no numbers $\equiv 0 \pmod{p}$,
it follows from Lemma \ref{cauchy_davenport_lemma} that  
$$
|L^\#_p|\ +\ |C^*_p|\ \leq\ p; 
$$
and so,
$$
{1 \over |L^\#_p|} + {1 \over |C^*_p|}\ \geq\ {1 \over |L^\#_p|} + {1 \over 
p-|L^\#_p|}\ \geq\ {4 \over p}.
$$
From this inequality we deduce
\begin{eqnarray}\label{LC_sum}
\sum_{p \leq Q} (\log p) \left ( {1 \over |L^\#_p|} + {1 \over |C^*_p|} \right )
&\geq& \sum_{p \leq Q} (\log p) \left ( {1 \over |L^\#_p|} + 
{1 \over p - |L^\#_p|} \right ) \nonumber \\
&\geq& \sum_{p \leq Q} {4 \log p \over p}\nonumber \\
&=& 2\log x + O(\log\log x).
\end{eqnarray}
Now, from Corollary \ref{gallagher_sieve2} with $J = L^\#$ and $J=C^*$ we deduce
$$
\sum_{p \leq Q} {\log p \over |L^\#_p|}\ <\ \log x\ +\ O(\log\log x)\ \ {\rm and\ \ }
\sum_{P \leq Q} {\log p \over |C^*_p|}\ <\ \log x\ +\ O(\log\log x).
$$
Combining these two upper bounds with (\ref{LC_sum}), we have that 
(\ref{LCQ_sum1}) is satisfied; and
$$
\sum_{p \leq Q} (\log p) \left ( {1 \over |L^\#_p|} + {1 \over p - |L^\#_p|}
\right )\ =\ 2\log x\ +\ O(\log\log x).
$$  
This equation and (\ref{LCQ_sum1}) together imply that
\begin{eqnarray}
\sum_{p \leq Q} {\log p \over p - |L^\#_p|} &=& 
\sum_{p \leq Q} (\log p) \left ( {1 \over |L^\#_p|} + {1 \over p - |L^\#_p|} \right )
\ -\ \sum_{p \leq Q} {\log p \over |L^\#_p|} \nonumber \\
&=& \log x\ +\ O(\log\log x), \nonumber
\end{eqnarray}
which gives that the first part of $(\ref{LCQ_sum2})$ is satisfied.
The second part of (\ref{LCQ_sum2}) is satisfied by applying the same
argument.
\bigskip

\section{Proof of Proposition \ref{variance_prop}.}
\bigskip

Let
$$
V_p(r) = \left ( J(r,p)\ -\ {|J| \over p-|K_p|} \right )^2. 
$$
The sum we wish to bound from above is as follows:
\begin{eqnarray}\label{UV_equation}
X &:=& \sum_{p \leq Q} (\log p) \sum_{r \in G_p} J(r,p)^2\nonumber \\
&\leq& \sum_{p \leq Q} (\log p)\sum_{r \in G_p} 
\left ( V_p(r)\ -\ 2{J(r,p)|J| \over p - |J_p|}\ +\ 
{|J|^2 \over (p - |J_p|)^2}\right ) \nonumber \\
&<& \sum_{p \leq Q} (\log p) \sum_{r \in J_p } V_p(r)
\ +\ |J|^2\sum_{p \leq Q} (\log p) \sum_{r \in G_p}
{1 \over (p - |J_p|)^2}. 
\end{eqnarray}

Now, we have that
$$
Y\ :=\ 
\sum_{p \leq Q} (\log p) \sum_{r \in \overline{(-K)_p}} V(r,p)\ =\ E_1\ -\ 2 E_2\ +\ E_3,
$$
where
\begin{eqnarray}
E_1 &=& \sum_{p \leq Q} (\log p) \sum_{r \in \overline{(-K)_p}} J(r,p)^2 
\ =\ \sum_{p \leq Q} (\log p) \sum_{r \in J_p} J(r,p)^2 \nonumber \\
&=& |J|^2 (\log x + O(1))\nonumber \\
E_2 &=& \sum_{p \leq Q} (\log p) \sum_{r \in \overline{(-K)_p}} 
{J(r,p)|J| \over p - |J_p|}\ =\ \sum_{p \leq Q} (\log p) 
\sum_{r \in J_p} {J(r,p) |J| \over p - |J_p|} \nonumber \\
&=& |J|^2 \sum_{p \leq Q} {\log p \over |p - K_p|} 
\ =\ |J|^2\left ( \log x\ +\ O(\log\log x) \right ) \nonumber \\
E_3 &=& |J|^2 \sum_{p \leq Q} {\log p \over p - |J_p|}\ =\ E_2. \nonumber
\end{eqnarray}
Note that the upper bound we derived for $E_1$ comes from Lemma \ref{gallagher_sieve1},
together with the fact that $\log Q = (\log x)/2 + O(\log\log x)$; 
the equation for $E_2$ comes from Lemma \ref{cauchy_davenport_consequence};
and, the switching of some of the above sums from a sum over $r \in \overline{(-K)_p}$
to $r \in J_p$ is justified since $J_p \subseteq \overline{(-K)_p}$,
by (\ref{L_subset_C}) and (\ref{C_subset_L}).

It follows that 
$$
Y\ =\ E_1 - E_2\ =\ O \left (|J|^2 \log\log x \right ).
$$

Substituting this into (\ref{UV_equation}) gives
\begin{eqnarray}
X &<& Y + |J|^2 \sum_{p \leq Q} (\log p) 
\sum_{r \in G_p} {1 \over (p-|J_p|)^2}\nonumber \\
&=& |J|^2 \left ( 
\sum_{p \leq Q} {(\log p) |G_p| \over (p-|J_p|)^2} + O(\log\log x)\right ),
\nonumber
\end{eqnarray}
which proves the Proposition.
\bigskip

\section{Proof of Lemma \ref{translation_invariant_lemma}.}

For an integer $h$, let $S(h)$ denote the symmetric difference between
$(J-hk)_p$ and $(J-(h-1)k)_p$.  We note that $|S(h)| = |S(0)|$.

Now, since
$$
J_p\ \setminus\ (J-hk)_p\ \ \subseteq\ \ S(h)\ \cup\ 
\bigl ( J_p\ \setminus\ (J-(h-1)k)_p \bigr ),
$$
it follows that 
$$
|J_p\ \setminus\ (J-hk)_p|\ \leq\ |S(0)|\ +\ |J_p\ \setminus\ (J-(h-1)k)_p|.
$$
For $h \geq 1$ a simple induction argument then shows that
\begin{equation}\label{setminus_bound}
|J_p\ \setminus\ (J-hk)_p|\ \leq\ h |S(0)|.
\end{equation}
Now, from (\ref{L_subset_C}) and (\ref{C_subset_L}) we deduce that
$J_p, (J+k)_p\ \subseteq\ \overline{(-K)_p}$, which gives: 
\begin{eqnarray}
|J_p\ \setminus\ (J+k)_p|\ \leq\ 
|\overline{(-K)_p}\ \setminus\ 
(J+k)_p | &=& |\overline{(-K)_p}|\ -\ |(J+k)_p|\nonumber \\
&=& p\ -\ |K_p|\ -\ |J_p|; \nonumber
\end{eqnarray}
and
$$
|(J+k)_p\ \setminus\ J_p|\ \leq\ |\overline{(-K)_p}\ \setminus\ 
J_p |\ =\ p\ -\ |K_p|\ -\ |J_p|.
$$
Thus,
\begin{equation}\label{setminus_bound2}
|S(0)|\ =\ |J_p\ \setminus\ (J+k)_p|\ +\ |(J+k)_p\ \setminus\ J_p|
\ \leq\ 2 ( p\ -\ |K_p|\ -\ |J_p|).
\end{equation}
From this and the fact that 
$$
p\ -\ |K_p|\ =\ |\overline{(-K)_p}|\ \leq\ |J_p|,
$$
we deduce  
\begin{eqnarray}
\sum_{p \leq Q} (\log p) {( p\ -\ |K_p|\ -\ |J_p| )
\over ( p\ -\ |K_p|)^2} &\leq& \sum_{p \leq Q} (\log p) 
{( p\ -\ |K_p|\ -\ |J_p|) \over |J_p| (p\ -\ |K_p|)} \nonumber \\
&=& \sum_{p \leq Q} (\log p) \left ( {1 \over |J_p|}\ -\ {1 \over 
p - |K_p|} \right )\nonumber \\
&=& O(\log\log x), \nonumber
\end{eqnarray}
by Lemma \ref{cauchy_davenport_consequence}.  
From this, (\ref{setminus_bound}), and (\ref{setminus_bound2}), we 
deduce 
\begin{eqnarray}
\sum_{p \leq Q} (\log p) { |J_p\ \setminus\ (J - jk)_p| \over 
(p\ -\ |K_p|)^2} &\leq& 2j \sum_{p \leq Q} (\log p) {(p\ -\ |K_p|\ -\ |K_p|)
\over (p - |K_p|)^2}\nonumber \\
&=& O(j\ \log\log x), \nonumber
\end{eqnarray}
which proves the Lemma.

\end{document}